\documentclass[12pt]{amsart}
\usepackage{amssymb}
\usepackage{graphicx}
\theoremstyle{plain}
\newtheorem{theorem}{Theorem}

\newcommand{\ga}{\alpha}

\newcommand{\cala}{\mathcal{A}}

\newcommand{\calg}{\mathcal{G}}

\newcommand{\calm}{\mathcal{M}}

\newcommand{\comp}{\mathrel{\scriptstyle\circ}}

\newcommand{\pd}[1]{\frac{\partial}{\partial{#1}}}

\newcommand{\Aut}{\operatorname{Aut}}

\newcommand{\C}{\mathbb{C}}

\newcommand{\coker}{\operatorname{coker}}

\newcommand{\Eb}{\mathbb{E}}

\newcommand{\Fl}{\operatorname{F\ell}}

\newcommand{\Fbb}{\mathbb{F}}

\newcommand{\g}{\mathfrak{g}}
\newcommand{\gl}{\mathfrak{gl}}
\newcommand{\GL}{\operatorname{GL}}

\newcommand{\Hbb}{\mathbb{H}}

\newcommand{\Hom}{\operatorname{Hom}}

\newcommand{\ind}{\operatorname{index}}
\newcommand{\inv}{^{-1}}

\newcommand{\Lie}{\mathcal{L}}

\newcommand{\liex}{\Lie_{X}}

\newcommand{\M}{\mathcal{M}}

\newcommand{\Ocal}{\mathcal{O}}

\newcommand{\Orbit}{\operatorname{Orbit}}

\newcommand{\Pont}{\operatorname{Pont}}

\newcommand{\pt}{\operatorname{pt}}

\newcommand{\R}{\mathbb{R}}

\newcommand{\SO}{\operatorname{SO}}

\newcommand{\Spin}{\operatorname{Spin}}

\newcommand{\SU}{\operatorname{SU}}

\newcommand{\frakt}{\mathfrak{t}}

\newcommand{\tr}{\operatorname{tr}}

\begin{document}
\title{The life and works of Raoul Bott}
\author{Loring W. Tu}
\date{\today}
\address{Department of Mathematics \\
     Tufts University\\
     Medford, MA 02155-7049}
\email{ltu@tufts.edu}
\maketitle

In a career spanning five decades, Raoul Bott has wrought profound 
changes on the landscape of geometry and topology.  It is a daunting 
task to improve upon his own reminiscences \cite{dioszeger}, \cite{auto}, 
\cite{recoll} and 
commentaries on papers \cite{collected}, 
punctuated as they are by insight, colorful 
turns of phrases and amusing anecdotes.  This article, based on his 
writings and on interviews and correspondences with people who have known 
him,
aims to 
serve as an introduction to his life and works so far.

\section{Early years}

Raoul Bott was born in Budapest in 1923.  His lineage fully 
reflects the geopolitical complexity of the region at the time.  
His mother's family was Hungarian and Jewish, while his father's side 
was Austrian and Catholic.  
His parents divorced soon after his birth, so he grew up with his mother 
and stepfather.
Raised as a Catholic, Raoul spent his childhood 
and adolescence in Slovakia, which seventy years later, 
after alternating between Hungary and 
Czechoslovakia, is 
today an independent country. 

In the first five years of school Raoul was not a good student.  This should 
give comfort to all parents of late bloomers.  In fact, he 
did not earn a single A except in singing and in German.  Nonetheless, 
he showed an early talent for breaking rules and for generating 
sparks---electrical sparks, that is, rigged up with wires, fuse 
boxes, vacuum tubes, and transformers.  The schools were formal and 
strict, and one could get slapped or have one's ears 
pulled for misbehaving.  For a budding original thinker, Raoul 
survived the schools relatively unscathed.  He recalls a friar hitting him on 
the hand once and a teacher cuffing his ear another time, for horsing 
around too much. 

It was by all accounts an idyllic existence, complete with a family 
villa, English governesses, and music lessons.  This 
world came to an abrupt halt in 1935, when his mother died of cancer.
In time his stepfather remarried.   

Raoul's experimental talent found its full flowering in adolescence.  
He and a kindred spirit Tomy Hornak built a small 
box with a slit for coins.  When someone dropped a coin through the slit, a 
display lit up saying ``Thank you.''  In this way they funded their early 
experiments.

Raoul struggled with some subjects in school and a tutor was hired to 
help him a few hours a week in his house.  At the time Raoul
 and Tomy had built a gadget
to communicate by Morse code.
As he was being tutored, he would 
hold the gadget under the table and Tomy would be sitting in the basement.
Raoul received the code by getting short and 
long electric shocks in his hand. 
He then responded by pressing a 
button to light up a bulb in the basement.  
While the tutor believed that his student was listening intently to 
the lesson, Raoul was chatting away in Morse code under the table.
In retrospect, Bott calls this his first attempt at e-mail.

\section{Canada}

In 1938, with Hitler's ascendancy and
Germany's march into Czechoslovakia,  Bott's stepparents flew him to 
the safety of England and enrolled him in an English boarding school.
Since they had only transit visas for England, the following 
year they headed for Canada, a country that to this day has been
extraordinarily welcoming to refugees and immigrants from 
around the world.

In the fall of 1941, after a rigorous year of preparatory studies in 
Ontario, Raoul Bott found himself at McGill University in Montreal.  Given 
his electrical know-how, he chose, not surprisingly, electrical engineering 
as his major.   His grades were respectable, but as he recalls in 
\cite{auto}, 
he was more interested 
in upholding the ``engineering tradition of hard drinking, loud, boisterous, 
mischievous, and macho behavior.''  Mathematics was his best subject; 
still, it was mathematics in the engineering sense, not the kind of pure 
reasoning for which he became so well known years later.

With his European flair, his 6 ft.\ 2 in.\ frame, and the conspicuous fur cap 
he often wore, Bott stood out from the crowd at McGill.  
When friends asked him where he was from, he said 
from Dioszeg, Czechoslovakia, and he added facetiously,
where he ``was a Count.''  After that, 
everyone called him the Count.

The Count sometimes spoke a very foreign tongue.
In the streetcars of Montreal, Raoul and his roommate Rodolfo
Gurdian would occasionally engage in 
a deliberately loud and animated conversation.  
Nothing they said made sense, for 
they were making up the language as they went along.  
From the corners of their eyes, they enjoyed watching the quizzical 
expressions on the faces of the surrounding passengers, who were trying hard 
to figure out what language the two of them were speaking.

Bott loved the opera, but as a penniless student how was he to afford it?  One
time the famous tenor Ezio Pinza came to sing in His Majesty's Theater, 
the opera house of Montreal in the Forties.  
For this occasion, Bott dressed up in his Sunday best and went to the theater.  
When the man at the entrance stopped him, Bott told him he couldn't do 
this because he was Ezio Pinza's nephew.  Bott said it with such assurance 
that the man let him in.  After that, Bott could go to all the shows at this 
theater for free.  

Bott's roommate Rodolfo, equally penniless, also loved the opera.  
But Rodolfo did not have the nerve to sneak into the theater.  
When the opera Carmen was playing, Rodolfo was very eager to attend.  
Bott magnanimously invited him.  By then, the ticket taker knew Bott 
very well, but he stopped Rodolfo at the entrance.  Bott turned around and 
intoned in his authoritative voice, ``It's all right.  He can come in.''  
Without any hesitation 
the ticket taker obeyed the order of this ``nephew'' of Ezio Pinza.

One New Year's Day, Raoul, Rodolfo, and some friends went to Mont 
Tremblant, a winter resort north of Montreal.  In the most prominent and 
expensive hotel, a big celebration was going on.  Somehow, to the envy of 
his friends, Raoul sneaked in.  A little later, Raoul was standing on the 
balcony, looking down contemptuously at his friends and
showing them a chicken leg he was eating.  After he finished it,
he threw the bone, with disdain, to his hungry friends.

(Old habits die hard.  In 1960 Bott, by then a full professor at 
Harvard, was in India with Michael Atiyah, both giving lectures as 
guests of the Tata Institute of Mathematics.  One day, as 
they walked in the streets of New Dehli, they passed by a big celebration.  
Bott decided 
to slip in uninvited, dragging Atiyah along with him. 
Atiyah, a professor at Oxford who was later 
anointed Sir Michael by the Queen and elected President of the Royal Society, 
was at first discomfited, but soon joined whole-heartedly in the 
festivities. 
They had a rousing time, sharing in the general merriment of complete strangers.)   

Upon graduation, Bott joined the army, but
the atomic bomb at Hiroshima put an end to his military career after only 
four months.  He 
entered a one-year Master's program in the Engineering Department at 
McGill.  Gradually it dawned on him that his interest lay more in 
mathematics than in engineering, and he produced a very mathematical 
master's thesis on ``impedance matching,'' which he said, ``the department 
accepted with some misgivings and about whose mathematical rigor I have 
doubts to this very day.''

At McGill Raoul met his future wife, Phyllis, an English literature major 
from the West Indies. 
Today, Phyllis remembers Raoul's first marriage proposal.  At the time he 
was doing his short stint in the army.  In full 
uniform, he said, ``Would you marry me?  Because if you do, the army will 
pay me more money.''  And then pointing through the window to his little 
room, he added, ``And we could be living there.''  The proposal was not 
accepted.  But two years later, they married.  The Botts have been together 
ever since, and now have four children and
eight grandchildren.  They celebrated their golden anniversary in 1997.

\section{Sermon}

While in the master's program in engineering at McGill, 
Bott floundered in trying to decide on the general direction of his career. 
Thirty years later, Bott was asked to deliver 
a sermon at Harvard's Memorial Chapel.  As he  
discussed the biblical passage of Eli, the wise man who counseled the 
young Samuel (1 Samuel 3: 3--6, 8--10), he reflected on the pivotal moment 
in his life that launched his mathematical career.  His description of his 
own Eli deserves to be read in the original:

\begin{quote}
And so when I saw the two readings we just heard juxtaposed in a 
Scripture Service, I could not resist them.  For they are appropriate to 
all of us, whether called to high causes or to lowly ones.  And they are 
maybe especially appropriate to the young people of today in their search 
of their destiny.

For surely there never has been a time when our young people have been 
given \emph{such freedom and therefore such responsibility} to find this 
destiny.

But how are we to know where we are called?  And how are we to know who is 
calling us?  These are questions beyond a mathematician's ken.  There are 
some who seem to have perfect pitch in these matters.  There are many more 
who might think that they have.  But with most of us, it is as it was with 
Samuel, and we are then truly blessed to have an advisor such as Eli.  He 
stands for all of us Teachers as an example.  For apart from communicating 
our call to our students, we should try and help them above all to discern 
theirs.

I well remember my Eli.  He was the Dean of the Medical School at McGill 
and I approached him for help in entering the medical school there, when in 
1945 the atomic bomb unexpectedly put an end to the war and to my four-month 
old career in the Canadian Infantry.

The Army very wisely decided to get rid of such green recruits as soon as 
possible, and so we all again found ourselves quite unexpectedly in charge 
of our own lives.  I had graduated in engineering earlier that year but had 
already decided against that career.

The Dean greeted me very cordially and assured me that there was a great 
need for technically trained doctors.  But, he said, seating me next to 
him, first tell me a little about yourself. Did you ever have any interest 
in botany, say, or biology? Well, not really, I had to admit.  How about 
chemistry --- Oh, I hated that course.  And so it went.  After a while he 
said, ``Well, is it maybe that you want to do good for humanity?''  And then, 
while I was coughing in embarrassment, he went on, ``Because they make the 
worst doctors.''

I thanked him, and as I walked out of his door I knew that I would start 
afresh and with God's grace try and become a mathematician.
\end{quote}

\section{Mathematical Career}

Initially Bott wanted to stay at McGill to do a mathematics Ph.D.
Because of his sketchy background, however, the McGill Math 
Department recommended that he pursue a Bachelor's degree in 
mathematics first.  It would have taken another three years.
Sensing his disappointment, Professor Williams of McGill then 
suggested Carnegie Tech (now Carnegie-Mellon University) 
to Bott, where John Synge was just forming a new graduate program and would 
need some students.

Synge received Bott warmly at Carnegie Tech, but as they read the rules of 
the program together, they found that Bott would have to spend three years 
taking courses in the newly minted master's program.  In a flash of 
inspiration, 
Synge said, ``Let's look at the Ph.D. program.''  It turned out to 
have hardly any requirements at all!  Normally the master's program is a 
prerequisite to the Ph.D. program, but perhaps recognizing a special gift 
in Bott, 
Synge put him in the Ph.D. program.  In just two years Bott would walk 
out with his degree.

Bott found the Carnegie Tech atmosphere exceedingly supportive.  The 
small coterie of mathematics students included Hans Weinberger, now 
at the University of Minnesota, and John Nash, an advanced 
undergraduate who after a thirty-year battle with schizophrenia
received the Nobel prize in 1994.  In later years Bott said of Carnegie 
Tech, ``Being a brand 
new graduate program, they hadn't learned yet how to put hurdles in front 
of graduate students.'' 
Bott considers himself very fortunate to have an advisor in R. J. Duffin, 
for Duffin treated him as an equal from the very outset and together they 
published two papers on the mathematics of electrical networks.

The first of these two papers, on impedance functions [1], 
so impressed Hermann Weyl that 
he invited Bott to the Institute for Advanced Study in 1949.  Thus began
Bott's initiation into the mysteries of algebraic topology.   
Apart from Weyl, among his 
main teachers were N.~Steenrod, E.~Specker, K.~Reidemeister, and M.~Morse. 
 Of Ernst Specker, 
Bott said in \cite{specker}, ``I bombarded Ernst with so many 
stupid questions that in 
desperation he finally imposed a fine of 25 cents on any conjecture he 
could disprove in less than five minutes. This should give you some idea of 
the inflation of the past thirty years and also help to explain Ernst's 
vast fortune at this time.''

At the time Norman Steenrod was writing his 
classic book on the topology of fiber bundles and teaching a course based 
on it.  This course greatly influenced 
Bott's mathematical development.  

Bott describes Steenrod with admiration as 
someone who treated high and low alike, with equal respect.  At Princeton, 
the graduate students could be intimidating, because they knew so much, and 
they let you know it.  Steenrod, on the other hand, was different.
In spite of his stature in the 
mathematical community, he put everyone at ease.  
In seminars Steenrod did not 
hesitate to ask the most basic questions.  This was quite often a boon to 
the others in the audience, too intimidated and too befuddled to ask the 
questions themselves.  

After two years at the Institute, Bott went to the University of 
Michigan.  In 1959 he became a professor at Harvard, where he has 
remained since.  In 1999 Bott formally retired from teaching.  He is now 
William Casper Graustein Research Professor at Harvard.

\section{Dunster House}

An unusual item 
in the curriculum vitae of Raoul Bott, for a mathematician at least, 
was his tenure as the Master of Dunster House in 1978--84.
At Harvard the undergraduates live in social units called ``Houses,'' 
modelled somewhat after the Colleges at Oxford and Cambridge.  A House is 
more than a place to sleep; it is a way to create a sense of a small community 
within a large university.  Each 
House has its own dining hall, dormitories, social activities, and a staff 
headed by a Master.  The academic staff consists of a bevy of resident and 
non-resident tutors.
    
Whether out of a lack of interest or a perceived mismatch of 
temperament, pure mathematicians are 
rarely called to be Masters of the undergraduate Houses.  In 1978, in a break 
with tradition, 
the President of Harvard University appointed
Bott the Master of Dunster House.  This entailed living in 
the Master's Residence in the midst of three hundred undergraduates.  
Bott's gregariousness was a good match for the post.

Every year the Houses compete in a water-raft race on the Charles River.  
This is no gentleman's canoe race as practiced in England.  In the 
Harvard version, attacks on other Houses' rafts are condoned, even 
encouraged.  One year the Lowell House team had its Master at the helm, 
resplendent in an admiral's hat. 
Bott, commanding the Dunster House 
armada, saw the beautiful hat.  He 
hollered, ``Get me that hat!''  Now, this is the sort of order 
undergraduates love to obey.  In no time the Dunster students 
had paddled to the Lowell raft.  A struggle ensued, and 
like any good pirates, the Dunster contingent captured the admiral's hat.  It was 
later hung, as a trophy, high in the ceiling of the Dunster House Dining Hall.

Showing true House spirit, 
the Dunster House Crew Team had its official team T-shirt 
emblazoned with  
``Dunster House,'' a pair of oars, 
and the exhortation: 
``Raoul, Raoul, Raoul your Bott.'' 

The Harvard Houses have counterparts at Yale, where they are called 
Colleges.  A friendly rivalry has always existed between these two august 
institutions, and it extends to the Houses and Colleges.  Some of the 
Houses at Harvard even have ``sister Colleges'' with which they are 
loosely affiliated.  They would, for example, visit each other during the 
Harvard-Yale football games.    

In the aftermath of the Sixties, many of the 
traditions at the Ivy League universities, such as the dress code and  
the parietal rules, have gone by the wayside, and for a number of 
years Dunster House had not had contact with Berkeley College, 
its sister College at Yale.  
One year the Berkeley College Master, a distinguished historian,
decided to revive the tradition.  He 
wrote to Bott suggesting a visit to Dunster House during the weekend of 
the Harvard-Yale football game.  Bott readily agreed, but decided to make 
the occasion a memorable one. 
Why not fool the Yalies into thinking that Harvard has kept up, at least to a 
certain point, the 
Oxbridge tradition of High Table and 
academic gowns at dinner? Why not show that, perhaps, 
Dunster House was more ``civilized'' 
than its Yale counterpart?  With enthusiasm, the Dunster House undergraduates 
all supported 
the idea.  

On the appointed day, the Dunster House Dining Hall was transformed from a 
cafeteria into a hallowed hall, complete with linen, waiters and waitresses, 
and even 
a wine steward wearing a large medal.  Unlike on a normal day, there were 
no T-shirts or cut-offs 
in sight.  Every tutor 
was attired in a black academic gown.  An orchestra sat in waiting.  When 
the Yale Master and his tutors arrived, 
Bott asked, with a straight face, ``Where are 
your gowns?''  Of course, they didn't have any.  ``Well, no problem, you 
could borrow some of ours.''  So the Dunster 
tutors led them to some gowns that had 
just been lent from Harvard's Chapel.  As Bott entered the Dining Hall with 
his guests, trumpets blared forth and the orchestra started playing.  The 
undergraduates were already seated, looking prim, proper, and serious.  
Bott and his tutors dined with the Yale visitors at a High Table, 
on a stage especially 
set up for this occasion.  The orchestra serenaded the diners with music.  
Everything went according to plan.  But the Yale Master, ever sharp, had 
the last laugh.  He opened his speech by saying, ``I'm glad to see that 
culture has finally migrated from New Haven to Harvard.''   

\section{Bott as a teacher}

Bott's lectures are legendary for their seeming ease of comprehension.  
His style is typically the 
antithesis of the Definition-Theorem-Proof approach so favored among 
mathematical speakers.  Usually he likes to discuss a simple key example 
that encapsulates the essence of the problem.  Often, as if by magic, a 
concrete formula with transparent significance appears.

At a reception for new graduate students at Harvard, he once shared his 
view of the process of writing a Ph.D. thesis.  He said it is like doing a 
homework problem; it's just a harder problem. You try to understand 
the problem thoroughly, from every conceivable angle.  Much of the thesis 
work is perseverance, as opposed to inspiration.  Above all, ``make the 
problem your own.''  

Many of his students testify to his warmth and humanity, but he also expects 
the students to meet an exacting standard.  He once banned the word 
``basically'' from an advisee's vocabulary, 
because that word to Bott signifies that some details 
are about to be swept under a rug.   

This insistence on thoroughness and clarity applies to his own work as well.
 I. M. Singer remarked that in their younger days, whenever they had a 
 mathematical discussion,
 the most common phrase Bott uttered was ``I 
 don't understand,'' and that a few months later Bott would emerge with a 
 beautiful paper on precisely the subject he had repeatedly not understood.

Seminar speakers at Harvard tend to address themselves to the experts in 
the audience.  But like 
Steenrod, Bott often interrupts the speakers with the most basic questions, 
with the salutary effect of slowing down the speakers
and making them more intelligible to lesser 
mortals. 

At Michigan and Harvard, Bott directed over 36 Ph.D. theses.  Some of his 
students have become luminaries in their own right: Stephen Smale 
and Daniel Quillen received the Fields Medal 
in 1966 and 1978 respectively, 
and Robert MacPherson the National Academy of Science 
Award in Mathematics in 1992.   The 
following is, I hope, the complete list 
of his Ph.D. students:

\medskip

\footnotesize
\begin{tabbing}
Year \quad\quad \= Alexandre, Antonio Franco \quad \=        \kill
{\bf Year} \quad\quad \> {\bf Ph.D.~Student}\quad\quad \> {\bf Dissertation Title} \\
\\
\rule{0pt}{15pt} 
1957 \> Smale, Stephen \> \parbox[t]{2.5in}{\raggedright 
Regular Curves on Riemannian Manifolds} \\
\rule{0pt}{15pt} 
1961 \> Edwards, Harold \> \parbox[t]{2.5in}{\raggedright Application of Intersection Theory to 
Boundary Value Problems} \\ 
\rule{0pt}{15pt} 
\>Curtis, Edward \> \parbox[t]{2.5in}{\raggedright The Lower Central Series for Free Group 
Complexes}\\
\rule{0pt}{15pt} 
1963 \> Conlon, Lawrence \> \parbox[t]{2.5in}{\raggedright 
Spaces of Paths on a Symmetric Space} \\ 
\rule{0pt}{15pt} 
\> Zilber, Joseph Abraham \> \parbox[t]{2.5in}{\raggedright 
Categories in Homotopy Theory}\\
\rule{0pt}{15pt}
1964 \> Holzsager, Richard Allan \>  \parbox[t]{2.5in}{\raggedright
Classification of Certain Types of Spaces}\\
\rule{0pt}{15pt}
\>Quillen, Daniel \> \parbox[t]{2.5in}{\raggedright Formal Properties of Over-Determined 
Systems Of Linear Partial Differential Equations}\\
\rule{0pt}{15pt}
1965 \> Landweber, Peter S. \>  \parbox[t]{2.5in}{\raggedright
K\"unneth Formulas for Bordism Theories}\\
\rule{0pt}{15pt}
\> Lazarov, Connor \> \parbox[t]{2.5in}{\raggedright
     Secondary Characteristic Classes in $K$-theory} \\
\rule{0pt}{15pt}
1969 \> Brooks, Morris William  \> \parbox[t]{2.5in}{\raggedright The Cohomology of the Complement
 of a Submanifold} \\
\rule{0pt}{15pt}
\> Brown, Richard Lawrence  \> \parbox[t]{2.5in}{\raggedright Cobordism Embeddings and Fibrations
 of Manifolds}\\
\rule{0pt}{15pt}
1970 \> Blass, Andreas R. \> Orderings of Ultrafilters\\
\rule{0pt}{15pt}
\> MacPherson, Robert D. \> \parbox[t]{2.5in}{\raggedright
 Singularities of Maps and Characteristic Classes}\\
\rule{0pt}{15pt}
1973 \> Miller, Edward Y.  \> \parbox[t]{2.5in}{\raggedright Local Isomorphisms of Riemannian
Hermitian and Combinatorial Manifolds}\\
\rule{0pt}{15pt}
1974 \> Garberson, John Dayton \> \parbox[t]{2.5in}{\raggedright
 The Cohomology of Certain Algebraic Varieties}\\
\rule{0pt}{15pt}
1975 \> Mostow, Mark  \> \parbox[t]{2.5in}{\raggedright Continuous Cohomology of Spaces 
with Two Topologies}\\
\rule{0pt}{15pt}
\> Perchik, James \>  \parbox[t]{2.5in}{\raggedright
Cohomology of Hamiltonian and Related Formal Vector Field Lie Algebras}\\
\rule{0pt}{15pt}
1976 \> Weiss, Richard Simon \>  \parbox[t]{2.5in}{\raggedright
Refined Chern Classes for Foliations}\\
\rule{0pt}{15pt}
1977 \> Brooks, Robert \> \parbox[t]{2.5in}{\raggedright On the Smooth Cohomology of Groups of
Diffeomorphisms} \\
\rule{0pt}{15pt}
1981 \> Hingston, Nancy \>  \parbox[t]{2.5in}{\raggedright
Equivariant Morse Theory and Closed Geodesics} \\
\rule{0pt}{15pt}
1982 \> Gunther, Nicholas \> \parbox[t]{2.5in}{\raggedright
 Hamiltonian Mechanics and Optimal Control} \\
\rule{0pt}{15pt}
\> Laquer, Turner Henry  \> \parbox[t]{2.5in}{\raggedright Homogeneous Connections and Yang-Mills
Theory on Homogeneous Spaces} \\
\rule{0pt}{15pt}
1984 \> Ticciati, Robin  \> \parbox[t]{2.5in}{\raggedright Singular 
Points in Moduli Spaces of 
Yang-Mills Fields} \\
\rule{0pt}{15pt}
1985 \> Forman, Robin  \> \parbox[t]{2.5in}{\raggedright Functional Determinants and Applications
  to Geometry}\\
\rule{0pt}{15pt}
1986 \> Corlette, Kevin  \> \parbox[t]{2.5in}{\raggedright Stability and Canonical Metrics 
in Infinite Dimensions}\\
\rule{0pt}{15pt}
1987 \> Block, Jonathan  \> \parbox[t]{2.5in}{\raggedright Excision in Cyclic Homology 
of Topological Algebras} \\
\rule{0pt}{15pt}
1989 \> Kocherlakota, Rama  \> \parbox[t]{2.5in}{\raggedright Integral Homology of Real Flag
  Manifolds and Loop Spaces of Symmetric Spaces} \\
\rule{0pt}{15pt}
\> Morelli, Robert  \> \parbox[t]{2.5in}{\raggedright  Hilbert's Third 
Problem and the $K$-Theory
  of Toric Varieties}\\
\rule{0pt}{15pt}
\> Bressler, Paul \>  \parbox[t]{2.5in}{\raggedright
Schubert Calculus in Generalized Cohomology} \\
\rule{0pt}{15pt}
1991 \> Grossberg, Michael  \> \parbox[t]{2.5in}{\raggedright Complete Integrability and
  Geometrically Induced Representations} \\
\rule{0pt}{15pt}
1992 \> Weinstein, Eric  \> \parbox[t]{2.5in}{\raggedright Extension of Self-Dual Yang-Mills
  Equations across the Eighth Dimension} \\
\rule{0pt}{15pt}
\> Szenes, Andras  \> \parbox[t]{2.5in}{\raggedright The Verlinde Formulas and Moduli Spaces of
  Vector Bundles} \\
\rule{0pt}{15pt}
1993 \> Tolman, Susan \>  \parbox[t]{2.5in}{\raggedright
Group Actions and Cohomology} \\
\rule{0pt}{15pt}
1994 \> Teleman, Constantine  \> \parbox[t]{2.5in}{\raggedright Lie Algebra Cohomology and 
the Fusion Rules}\\
\rule{0pt}{15pt}
\> Costes, Constantine  \> \parbox[t]{2.5in}{\raggedright Some Explicit Cocycles for 
Cohomology Classes Of Group Diffeomorphisms Preserving a $G$-Structure}\\
\rule{0pt}{15pt}
2000 \> Bernhard, James  \> \parbox[t]{2.5in}{\raggedright Equivariant de Rham Theory 
and Stationary Phase Expansions}
\end{tabbing}

\normalsize
\medskip

\section{Honors and Awards}

Throughout his career, Bott has been showered with honors, awards, and 
prizes.  The more noteworthy awards include: Sloan Fellowship (1956--60),
 Veblen Prize of the American Mathematical Society (1964), 
Guggenheim Fellowship (1976), National 
Medal of Science (1987), Steele Career Prize of the American Mathematical 
Society (1990), and the 
Wolf Prize in Mathematics (2000).

He was twice invited to address the International Congress of 
Mathematicians, in Edinburgh in 1958 and in Nice in 1970. 

He was elected Vice-President of the American Mathematical Society in 
1974--75, Honorary Member of the London Mathematical Society (1976), 
Honorary Fellow of St.\ Catherine's College, Oxford (1985), and Honorary 
Member of the Moscow Mathematical Society (1997).  
He has been a member of the 
National Academy of Science since 1964 and the French Academy of Sciences 
since 1995.

In 1987 he gave the Convocation Address at McGill University.    
He has also received Honorary Degrees of Doctor of Science from the 
University of Notre Dame (1980), McGill University (1987), Carnegie Mellon 
University (1989), and the University of Leicester, England (1995).  

\bigskip
\noindent
{\bf Mathematical Works}
\medskip

The bibliography in Raoul Bott's Collected Papers \cite{collected} 
lists his publications, 
with some omissions, up to 1990.  At the end of this article, we complete 
that bibliography by listing the missing publications up to 1990 and the 
publications since then. 

When asked to single out the top three in the manner of an Olympic 
contest, he replied, ``Can I squeeze in another one?''  But after listing four 
as the tops, he sighed and said, ``This is like being asked to single out 
the favorites among one's children.''  In the end he came up with a top-five 
list, in chronological order:
\begin{enumerate}
\item[{[15]}] Homogeneous vector bundles,
\item[{[24]}] The periodicity theorem,
\item[{[51]}] Topological obstruction to integrability,
\item[{[81]}] Yang-Mills equations over Riemann surfaces, 
\item[{[82]}] The localization theorem in equivariant cohomology.
\end{enumerate}

To discuss only these five would not do justice to the range of his output.
On the other hand, it is evidently not possible to 
discuss every item in his ever-expanding opus.  As a compromise, I
asked him to make a longer list of all his favorite papers, without trying 
to rank them.  What follows is a leisurely romp 
through the nineteen papers he chose.  
My goal is to explain, as simply as possible, the main 
achievement of his own favorite papers.  
For this reason, the theorems, if stated at all, are often not 
in their greatest generality.   

\section{Impedance}

The subject of Raoul Bott's first paper [1] dates back to his engineering days.  
An
electrical network determines an impedance function $Z(s)$, which describes the 
frequency response of the network.  This impedance function $Z(s)$ is a 
rational function of a complex variable $s$ and is {\it positive-real} 
(p.r.) in 
the sense that it maps the right half-plane into itself.  An old question 
in electrical engineering asks whether conversely, given a positive-real 
rational function $Z(s)$, it is possible to build a network with $Z(s)$ as 
its impedance function.  In some sense O. Brune had solved this problem in 
1931, but Brune's solution assumes the 
existence of an ``ideal transformer,'' which in practice would have to be 
the size of, say, the Harvard Science Center.  The assumption of an ideal 
transformer renders Brune's algorithm not so practical, and it was Raoul's 
dream at McGill to remove the ideal transformer from the solution.

At his first meeting with his advisor Richard Duffin at Carnegie Tech, he 
blurted out the problem right away.  Many days later, after 
a particularly fruitless and strenuous discussion, Raoul went home and 
realized how to do it.  He called Duffin.  The phone was busy.  As it 
turned out, Duffin was 
calling him with exactly the same idea!  
They wrote up the solution to the long-standing problem in a joint paper, 
which amazingly took up only two pages.

\section{Morse theory}

As mentioned earlier, the paper on impedance so impressed Hermann Weyl 
that he invited Bott to the Institute for Advanced Study at Princeton in 1949.  
There Bott came into 
contact with Marston Morse.  Morse's theory of critical points would play 
a decisive role throughout Bott's career, notably in his work on homogeneous 
spaces, the Lefschetz hyperplane theorem, the periodicity theorem, 
and the Yang-Mills functional on a moduli space. 

In the Twenties Morse had initiated the study of the critical points of 
a function on a space and its relation to the topology of the space.
A smooth function $f$ on a smooth manifold $M$ has a \emph{critical point} 
at $p$ in $M$ if there is a coordinate system $(x_1, \ldots, x_n)$ at $p$ 
such that all the partial derivatives of $f$ vanish at $p$: 
\[
\dfrac{\partial f}{\partial x_i} (p) =0  \quad\quad \text{for all } i= 
1,\ldots, n.
\]
Such a critical point is \emph{nondegenerate} if the matrix of second 
partials, called the \emph{Hessian} of $f$ at $p$,
\[
H_p f= \left[ \dfrac{\partial ^2 f}{\partial x_i \partial x_j}(p) \right],
\]
is nonsingular.
The \emph{index} $\lambda (p)$ of a nondegenerate critical point $p$ is the 
number of negative eigenvalues of the Hessian $H_p f$; it is the number of 
independent directions along which $f$ will decrease from $p$.

If a smooth function has only nondegenerate critical points, we call it a 
\emph{Morse function}. 
The behavior of the critical points of a Morse function can be summarized 
in its \emph{Morse polynomial}:
\[
\M _t(f) := \sum t^{\lambda (p)},
\]
where the sum runs over all critical points $p$.

\begin{figure}[!h]
\begin{center}
\includegraphics[scale=0.25]{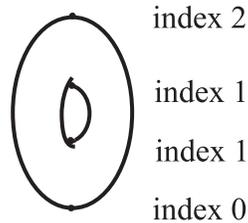}
 \caption{Critical points of the height function}
 \label{vertorus}
\end{center}
\end{figure}

A typical example of a Morse function is 
the height function $f$ of a torus standing 
vertically on a table top (Figure~\ref{vertorus}).

\medskip

\noindent
The height function on this torus has four critical points of index 0, 
1, 1, 2 respectively.  Its Morse polynomial is 
\[
\M_t (f) = 1 + 2t+t^2.
\]

For a Morse function $f$ on a compact manifold $M$, the fundamental results of 
Morse theory hinge on the fact 
that $M$ has the homotopy type of a CW complex with 
one cell of dimension $\lambda$ for each critical point of $f$ of index 
$\lambda$.
This realization came about in the early Fifties, due to the work of 
Pitcher, Thom, and Bott.

Two consequences follow immediately:
\begin{enumerate}
\item[i)] The weak Morse inequalities:
\[
\#\ \text{critical points of index } i \ge i\text{-th Betti number}.
\]
If 
\[
P_t (M) = \sum \dim H_i(M) t^i
\]
is the Poincar\'e polynomial of $M$, the Morse inequalities can be restated 
in the form 
\[
\M _t(f) \ge P_t(M),
\]
meaning that their difference $\M _t(f) - P_t(M)$  
is a polynomial with nonnegative coefficients.  This inequality provides a 
topological constraint on analysis, for it says that the $i$-th Betti 
number of the manifold sets a lower bound on the number of critical points 
of index $i$ that the function $f$ must have.
\item[ii)] The lacunary principle:  If no two 
critical points of the Morse function $f$ have consecutive indices, then 
\begin{equation} \label{perfect}
\M _t(f) = P_t(M).
\end{equation}
The explanation is simple: since in the CW complex of $M$ there are no two 
cells of consecutive dimensions, the boundary operator is automatically 
zero.  Therefore, the cellular chain complex is its own homology.
\end{enumerate}

A Morse function $f$ on $M$ satisfying \eqref{perfect} is said to be 
\emph{perfect}.  The height function on the torus above is a perfect Morse 
function.

\begin{figure}[!h]
\begin{center}
\includegraphics[scale=0.5]{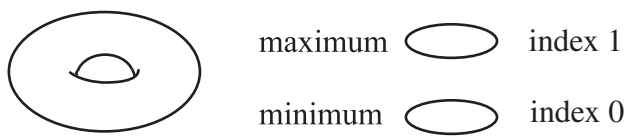}
 \caption{Critical manifolds of the height function}
 \label{hortorus}
\end{center}
\end{figure}

Classical Morse theory deals only with functions all of 
whose critical points are nondegenerate; in particular, the critical 
points must all be isolated points.  In many situations, however, the 
critical points form submanifolds of $M$.  For example, if the torus now 
sits flat on the table, as a donut usually would, then the height function 
has the top and bottom circles as critical manifolds 
(Figure~\ref{hortorus}).

\medskip

One of Bott's first insights was to see how to extend Morse theory to this 
situation.  In [9] he introduced the notion of a nondegenerate critical 
manifold:  a critical manifold $N$ is \emph{nondegenerate} if at any point 
$p$ in $N$ the Hessian of $f$ restricted to the normal space to $N$ is 
nonsingular.  The index $\lambda (N)$ of the nondegenerate critical 
manifold $N$ is then defined to be the number of negative eigenvalues of 
this normal Hessian; it represents the number of independent normal 
directions along which $f$ is decreasing.
For simplicity, assume that the normal bundles of the nondegenerate 
critical manifolds are all orientable.  
To form the Morse polynomial of 
$f$, each critical manifold $N$ is counted with its Poincar\'e polynomial; 
thus, 
\[
\M _t(f) := \sum P_t(N) t^{\lambda (N)},
\]
summed over all critical manifolds.  

With this definition of the Morse polynomial, Bott proved in [9] that if a 
smooth function $f$ on a smooth manifold $M$ has only nondegenerate 
critical manifolds, then the Morse inequality again holds:
\[
\M _t(f) \ge P_t(M).
\]

\section{Lie groups and homogeneous spaces}

In the Fifties Bott applied Morse theory with great success to the topology 
of Lie groups and homogeneous spaces.  
In [8] he showed how the diagram of a compact semisimple 
connected and simply connected group $G$ determines the integral homology 
of both the loop space $\Omega G$ and the flag manifold $G/T$, where $T$ 
is a maximal torus.

Indeed, Morse theory gives a beautiful CW cell structure on $G/T$, up to 
homotopy equivalence.  To 
explain this, recall that the adjoint action of the group $G$ on its Lie 
algebra $\g$ restricts to an action of the maximal torus $T$ on $\g$.  As 
a representation of the torus $T$, the Lie algebra $\g$ decomposes into a 
direct sum of irreducible representations
\[
\g = \frakt \oplus \sum E_{\alpha} ,
\]
where $\frakt$ is the Lie algebra of $T$ and each $E_{\ga}$ is a 
$2$-dimensional space on which $T$ acts as a rotation $e^{2\pi i \alpha 
(x)}$, corresponding to the root $\alpha (x)$ on $\frakt$.
The \emph{diagram} of $G$ is the family of parallel hyperplanes in $\frakt$ 
where some root is integral.  A hyperplane that is the zero set of a root 
is called a \emph{root plane}.

\begin{figure}[!h]
\begin{center}
\includegraphics[scale=0.4]{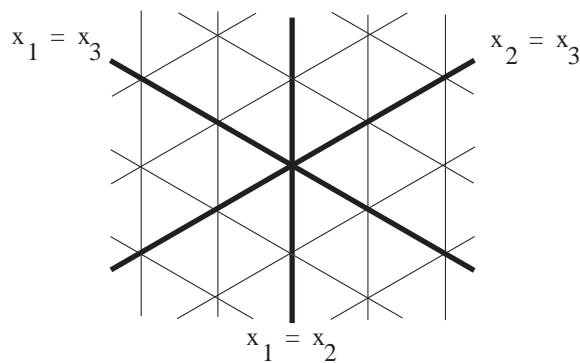}
 \caption{The diagram of SU(3)}
 \label{diagram}
\end{center}
\end{figure}

For example, for the group $G = \SU (3)$ and maximal torus
\[ 
T=\left\{ \left. 
\begin{bmatrix}
e^{2\pi i x_1} & & \\
& e^{2\pi i x_2} & \\
& & e^{2\pi i x_3}
\end{bmatrix}
\
\right| \ 
x_1+x_2+x_3 = 0, x_i \in \R \right\} .
\]
the roots are $\pm(x_1 - x_2), \pm(x_1 - x_3) , \pm(x_2 - x_3)$, 
and the diagram is the collection of lines in the plane $x_1+x_2+x_3=0$ in 
$\R^3$ as in Figure \ref{diagram}.  In this figure, the root planes are 
the thickened lines.

For $G=\SU (2)$ and
\[
T= \left\{ \left. 
\begin{bmatrix}
e^{2\pi i x} & 0 \cr
0 & e^{-2\pi i x}
\end{bmatrix}
\ \right| \ 
x \in\R \right\},
\]
the Lie algebra $\frakt$ is $\R$, the roots are $ \pm 2x$, and the adjoint 
representation of of $G$ on $\g = \R^3$ corresponds to rotations.  The root 
plane is the origin.

\begin{figure}[!h]
\begin{center}
\includegraphics[scale=0.4]{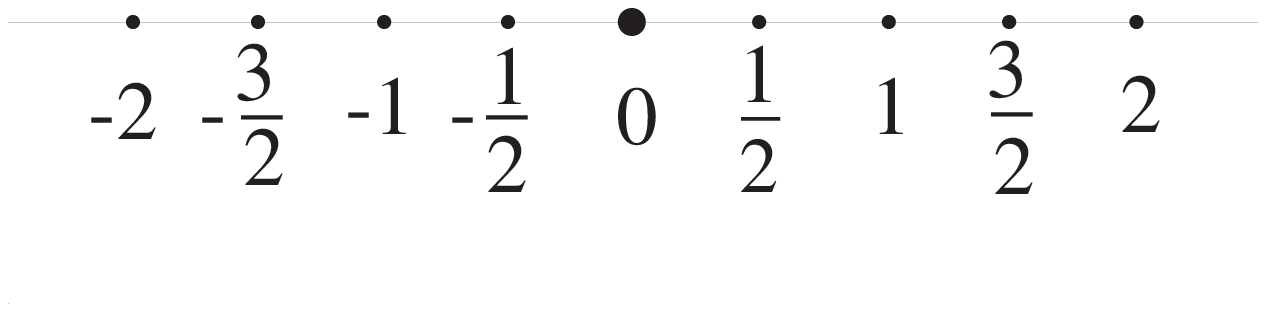}
 \caption{The diagram of SU(2)}
 \label{su2}
\end{center}
\end{figure}

A point $B$ in $\frakt$ is \emph{regular} if its normalizer has minimal 
possible dimension, or equivalently, if its normalizer is $T$.  It is well 
known that a point $B$ in $\frakt$ is regular if and only if it does not lie 
on any of the hyperplanes of the diagram.  If $B$ is regular, then the 
stabilizer of $B$ under the adjoint action of $G$ is $T$ and so the orbit 
through $B$ is $G/T$.

Choose another regular point $A$ in $\frakt$ and define the function $f$ 
on $\Orbit (B) = G/T$ to be the distance from $A$; here the distance is 
measured with respect to the Killing form on $\g$.  
Let $\{ B_i\}$ be all 
the points in $\frakt$ obtained from $B$ by reflecting about the root 
planes.  Then Bott's theorem asserts that $f$ is a Morse function on $G/T$ 
whose critical points are precisely all the $B_i$'s.  Moreover, the index 
of a critical point $B_i$ is twice the number of times that the line 
segment from $A$ to $B_i$ intersects the root planes.  This cell 
decomposition of Morse theory fits in with the more group-theoretic Bruhat 
decomposition.

\begin{figure}[!h]
\begin{center}
\includegraphics[scale=0.4]{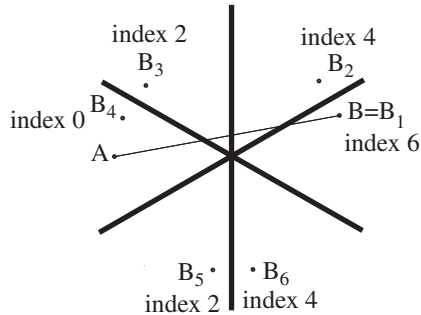}
 \caption{The flag manifold $\Fl (1,2,3)$}
 \label{flag}
\end{center}
\end{figure}

For $G=\SU (3)$ and  $T$ the set of diagonal matrices in $\SU(3)$, the 
orbit $G/T$ is the complex flag manifold $\Fl ( 1,2,3)$, consisting of all 
flags
\[
V_1 \subset V_2 \subset \C^3, \quad \dim_{\C} V_i = i.
\]
Bott's recipe gives 6 critical points of index $0,2,2,4,4,6$ respectively 
on $G/T$ (See Figure \ref{flag}).  By the lacunary principle, the Morse function 
$f$ is perfect.  Hence, the flag manifold $\Fl (1,2,3)$ has the homotopy 
type of a CW complex with one $0$-cell, two $2$-cells, two $4$-cells, and one $6$-cell.  
Its Poincar\'e polynomial is therefore
\[
P_t(\Fl (1,2,3)) = 1+2t^2 +2t^4 +t^6.
\]    

\section{Index of a closed geodesic}

For two points $p$ and $q$ on a Riemannian manifold $M$, the space 
$\Omega_{p,q}(M)$ of all paths from $p$ to $q$ on $M$ is not a 
finite-dimensional manifold.  
Nonetheless, Morse theory applies to this situation also, with a Morse 
function on the path space $\Omega_{p,q}$ given by the energy of a path:
\[
E(\mu )=\int_a^b \langle \frac{d\mu}{dt}, \frac{d\mu}{dt} \rangle dt.
\]
The first result of this infinite-dimensional Morse theory asserts that the 
critical points of the energy function are precisely the geodesics from 
$p$ to $q$.

Two points $p$ and $q$ on a geodesic are \emph{conjugate} if keeping $p$ 
and $q$ fixed, one can vary the geodesic from $p$ to $q$
through a family of geodesics. 
For example, two antipodal points on an $n$-sphere are conjugate points.  
The \emph{multiplicity} of $q$ as a conjugate point of $p$ is the 
dimension of the family of geodesics from $p$ to $q$.  On the $n$-sphere 
$S^n$, the multiplicity of the south pole as a conjugate point of the north 
pole is therefore $n-1$.

If $p$ and $q$ are not conjugate along the geodesic, then the geodesic is 
nondegenerate as a critical point of the energy function on $\Omega_{p,q}$. 
Its index, according to the Morse index theorem, is the number of 
conjugate points from $p$ to $q$ counted with multiplicities.

\begin{figure}[!h]
\begin{center}
\includegraphics[scale=0.4]{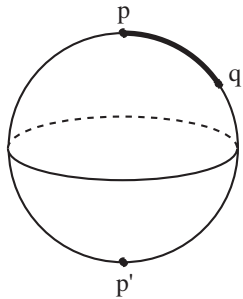}
 \caption{Geodesics on a sphere}
 \label{sphere}
\end{center}
\end{figure}

On the $n$-sphere let $p$ and $p'$ be antipodal points and $q \ne p'$.  The 
geodesics from $p$ to $q$ are $pq, pp'q, pqp'pq, pp'qpp'q, \dots$, of 
index $0, n-1, 2(n-1), 3(n-1), \dots$, respectively.  By the Morse index 
theorem the energy function on the path space $\Omega_{p,q}(S^n)$ has one 
critical point each of index $0, n-1, 2(n-1), 3(n-1), \dots$.  It then 
follows from Morse theory that $\Omega_{p,q}(S^n)$ has the homotopy type 
of a CW complex with one cell in each of the dimensions  
$0, n-1, 2(n-1), 3(n-1), \dots$.

Now consider the space $\Omega M$ of all smooth loops in $M$, that is, smooth 
functions $\mu: S^1 \to M$.  The critical points of the energy function 
on $\Omega M$ are again the geodesics, but these are now closed geodesics.  
A closed geodesic is never isolated as a critical point, since for any 
rotation $r: S^1 \to S^1$ of the circle, $\mu \comp r: S^1 \to M$ is still 
a geodesic.  In this way, any closed geodesic gives rise to a circle of 
closed geodesics.  When the Riemannian metric on $M$ is generic, the 
critical manifolds of the energy function on the loop space $\Omega$ will 
all be circles.  

Morse had shown that the index of a geodesic is the 
number of negative eigenvalues of a Sturm differential equation, 
a boundary-value problem of the form $Ly=\lambda y$, where $L$ 
is a self-adjoint second-order differential operator.
For certain boundary conditions, Morse had expressed the index in terms of 
conjugate points, but this procedure
does not apply to closed geodesics, which 
correspond to a Sturm problem with periodic boundary conditions.

In [14] Bott found an algorithm to compute the 
index of a closed geodesic.  He was then able to determine the behavior of the 
index when the closed geodesic is iterated.  Bott's method is in fact applicable 
to all Sturm differential equations.  And so in his paper he also gave a 
geometric formulation and new proofs of the Sturm-Morse separation, 
comparison, and oscillation theorems, all based on the principle that the 
intersection number of two cycles of complementary dimensions is zero if 
one of the cycles is homologous to zero.

\section{Homogeneous vector bundles}

Let $G$ be a connected complex semisimple Lie group, and $P$ a parabolic 
subgroup.  Then $G$ is a principal $P$-bundle over the homogeneous 
manifold $X= G/P$.
Any holomorphic representation $\phi : P \to \Aut (E)$ on a complex 
vector space $E$ induces a holomorphic vector bundle $\Eb$ over $X$:
\[
\Eb:= G \times_{\phi} E := (G \times E)/ \sim,
\]
where $(gp,e)\sim (g, \phi(p) e)$.  Then $\Eb$ is a 
holomorphic vector bundle over $X= G/P$.
A vector bundle over $X$
arising in this way is called a \emph{homogeneous vector 
bundle}.  
Let $\Ocal (\Eb)$ be the 
corresponding sheaf of holomorphic sections.
The homogeneous vector bundle $\Eb$ inherits a left $G$-action from the left 
multiplication in $G$:
\[
h.(g,e)=(hg, e)  \quad\quad \text{for } h,g \in G, e \in E.
\]
Thus, all the cohomology groups $H^q(X, \Ocal(\Eb))$ become $G$-modules.

In [15] Bott proved that if the representation $\phi$ is irreducible, the 
cohomology groups $H^q(X, \Ocal(\Eb))$ all vanish except possibly in one 
single dimension.  Moreover, in the nonvanishing dimension $q$, 
$H^q(X, \Ocal(\Eb))$ is an irreducible representation of $G$ whose highest 
weight is related to 
$\phi$.

This theorem generalizes an earlier theorem of Borel and Weil, who proved 
it for a positive line bundle.

In Bott's paper one finds a precise way of determining the nonvanishing 
dimension in terms of the roots and weights of $G$ and $P$.  Thus, on the 
one hand, Bott's theorem gives a geometric realization of induced 
representations, and on the other hand, it provides an extremely useful 
vanishing criterion for the cohomology of homogeneous vector bundles.

\section{The periodicity theorem}

Homotopy groups are notoriously difficult to compute.  For a simple space 
like the $n$-sphere, already, the higher homotopy groups exhibits no 
discernible patterns.  It was therefore a complete surprise in 1957, when 
Raoul Bott computed the stable homotopy groups of the classical groups and 
found a simple periodic pattern for each of the classical groups [24].

We first explain what is meant by the \emph{stable} homotopy group.
Consider the unitary group $U(n+1)$.  It acts transitively on the unit 
sphere $S^{2n+1}$ in $\C^{n+1}$, with stabilizer $U(n)$ at the point 
$(1,0,\dots,0)$. In this way, the sphere $S^{2n+1}$ can be identified with 
the homogeneous space $U(n+1)/U(n)$, and there is a fibering $U(n+1) \to 
S^{2n+1}$ with fiber $U(n)$.
By the homotopy exact sequence of a fibering, the following sequence is 
exact:
\[
\dots \to \pi_{k+1}(S^{2n+1}) \to \pi_k (U(n)) \to \pi_k (U(n+1)) \to
\pi_k (S^{2n+1}) \to \dots.
\]
Since $\pi_k (S^m) = 0$ for $m > k$, it follows immediately that as $n$ goes 
to infinity (in fact for all $n > k/2$), the $k$th homotopy group of the 
unitary group stabilizes:
\[
\pi_k (U(n)) = \pi_k (U(n+1)) = \pi_k(U(n+2)) = \dots.
\]
This common value is called the \emph{$k$th stable homotopy group} of the 
unitary group, denoted $\pi_k (U)$.

In the original proof of the periodicity theorem [24], Bott showed that in the 
loop space of the special unitary group $\SU(2n)$, the manifold of minimal 
geodesics is the complex Grassmannian 
\[
G(n,2n) =\frac{U(2n)}{U(n)\times U(n)}.
\]
By Morse theory, the loop space $\Omega \SU (2n)$ has the homotopy type of 
a CW complex obtained from the Grassmannian $G(n,2n)$ by attaching
cells of dimension $\ge 2n+2$:
\[
\Omega \SU (2n) \sim G(n,2n) \cup e_{\lambda} \cup \dots, 
\quad\dim e_{\lambda} \ge 
2n+2.
\]
It follows that
\[
\pi_k (\Omega \SU (2n)) = \pi_k (G(n,2n))
\]
for $n >> k$.

It is easily shown that
\[
\pi_k (\Omega \SU (2n)) = \pi_{k+1}(\SU(2n)) = \pi_{k+1}(U(2n)).
\]
Using the homotopy exact sequence of the fibering
\[
U(n) \to U(2n)/U(n)  \to G(n,2n),
\]
one gets
\[
\pi_k (G(n,2n)) = \pi_{k-1} (U(n)).
\]
Putting all this together, for $n$ large relative to $k$, we get
\[
\pi_{k-1} (U(n)) = \pi_k (G(n,2n)) = \pi_k(\Omega\SU(2n)) = \pi_{k+1}(U(2n)).
\]
Thus, the stable homotopy group of the unitary group is periodic of period 
$2$:
\[
\pi_{k-1} (U) = \pi_{k+1} (U).
\]

Applying the same method to the orthogonal group and the symplectic group, 
Bott showed that their stable homotopy groups are periodic of period $8$.

\section{Clifford algebras}

The \emph{Clifford algebra} $C_k$ is the algebra over $\R$ with $k$ 
generators $e_1, \dots, e_k$ and relations
\begin{align*}
e_i^2 &= -1   \quad \quad \text{for } i=1,\dots, k, \cr
e_ie_j &= - e_je_i \quad \quad \text{for all } i \ne j.
\end{align*}

The first few Clifford algebras are easy to describe
\[
C_0 = \R,\quad C_1 = \C,\quad C_2 = \Hbb =\{ \text{quaternions} \}.
\]
If $\Fbb$ is a field, denote by $\Fbb (n)$ the algebra of all $n\times n$ 
matrices with entries in $\Fbb$.  We call $\Fbb (n)$ a full matrix algebra.  It 
turns out that the Clifford algebras are all full matrix algebras or the 
direct sums of two full matrix algebras:
\[
\begin{array}{ccc|ccc|cc}
k & C_k && k & C_k && k & C_k \cr
\hline
0 & \R && 8 & \R (16) && 16 & \R(2^8) \cr
1 & \C && 9 & \C (16) && 17 & \C(2^8) \cr
2 & \Hbb && 10 & \Hbb (16) && 18 &  \vdots \cr
3 & \Hbb \oplus \Hbb && 11 & \Hbb(16) \oplus \Hbb(16) &&& \cr
4 & \Hbb (2) &&12 & \Hbb(32) &&& \cr
5 & \C (4) && 13  & \C(64)   &&& \cr
6 & \R (8) && 14  & \R(128) &&& \cr
7 & \R(8)\oplus\R(8) && 15 & \R(128)\oplus \R(128) &&& 
\end{array}
\]

This table exhibits clearly a periodic pattern of period $8$, 
except for the dimension increase after each period.
The $8$-fold periodicity of the Clifford algebras, long known to 
algebraists, is reminiscent of the $8$-fold 
periodicity of the stable homotopy groups of the orthogonal group.

In the early Sixties Michael Atiyah, Raoul Bott, and Arnold Shapiro 
found an explanation 
for this tantalizing connection. The link is provided by
a class of linear differential 
operators called the \emph{Dirac operators}.  The link between differential 
equations and homotopy groups first came about as a result of the realization 
that ellipticity of a differential operator 
can be defined in terms of the symbol of the differential operator. 

Suppose we can find $k$ real matrices $e_1, \dots, e_k$ of size $n \times 
n$ 
satisfying 
\[
e_i^2 = -1, \quad e_ie_j=-e_je_i \quad \text{for } i \ne j.
\]
This corresponds to a real representation of the Clifford algebra $C_k$.
The associated Dirac operator $D=D_{k,n}$ is the linear first-order 
differential operator
\[
D= I \pd{x_0} + e_1 \pd{x_1} + \dots + e_k \pd{x_k},
\]
where $I$ is the $n \times n$ identity matrix.  Such a differential 
operator on $\R^{k+1}$ has a \emph{symbol} $\sigma_D(\xi)$ obtained by 
replacing $\partial/\partial{x_i}$ by a variable $\xi_i$:
\[
\sigma_D(\xi) = I \xi_0 + e_1 \xi_1 + \dots + e_k \xi_k.
\]
The Dirac operator $D$ is readily shown to be \emph{elliptic}; this means 
its symbol
$\sigma_D(\xi)$ is nonsingular for all $\xi \ne 0$ in $\R^{k+1}$.
Therefore, when restricted to the unit sphere in $\R^{k+1}$, the symbol of 
the Dirac operator gives a map
\[
\sigma_D(\xi): S^k \to \GL (n, \R).
\]
Since $\GL (n,\R)$ has the homotopy type of $O(n)$, this map given by the 
symbol of the Dirac operator defines an element of the homotopy group
$\pi_k(\GL (n,\R)) = \pi_k (O(n))$.

The paper [33] shows that the minimal-dimensional representations of the 
Clifford algebras give rise to Dirac operators whose symbols generate the 
stable homotopy groups of the orthogonal group.  In this way, the $8$-fold 
periodicity of the Clifford algebras reappears as the $8$-fold periodicity 
of the stable homotopy groups of the orthogonal group.

\section{The index theorem for homogeneous differential operators}

The Sixties was a time of great ferment in topology and one of its 
crowning glories was the Atiyah-Singer index theorem.  Independently of 
Atiyah and Singer's work, Bott's paper [37] on homogeneous differential 
operators 
analyzes an interesting example where the analytical difficulties can be 
avoided by representation theory.

Suppose $G$ is a compact connected Lie group and $H$ a closed connected 
subgroups.  As in our earlier discussion of homogeneous vector 
bundles, a representation $\rho$ of $H$ gives rise to a vector bundle $G 
\times _{\rho} H$ over the homogeneous space $X=G/H$.  Now suppose $E$ 
and $F$ are two vector bundles over $G/H$ arising from representations of 
$H$.  Since $G$ acts on the left on both $E$ and $F$, it also acts on their 
spaces of sections, $\Gamma (E)$ and $\Gamma (F)$.  We say that a 
differential operator $D: \Gamma (E) \to \Gamma (F)$ is \emph{homogeneous} 
if it commutes with the actions of $G$ on $\Gamma (E)$ and $\Gamma (F)$.
If $D$ is elliptic, then its index
\[
\ind (D) = \dim \ker D - \dim \coker D
\]
is defined.

Atiyah and Singer had given a formula for the index of an elliptic operator 
on a manifold in terms of the topological data of the situation:  the 
characteristic classes of $E$, $F$, the tangent bundle of the base 
manifold, and the symbol of the operator.  In [37] Raoul Bott verified the 
Atiyah-Singer index theorem for a homogeneous operator by introducing a 
refined index, which is not a number, but a character of the group $G$.  
The usual index may be obtained from the refined index by evaluating at the 
identity.  A similar theorem in the infinite-dimensional case has recently 
been proven in the context of physics-inspired mathematics.

\section{Nevanlinna theory and the Bott-Chern classes}

Nevanlinna theory deals with the following type of questions:  Let $f:\C 
\to \C P^1$ be a holomorphic map.  Given $a$ in $\C P^1$, what is 
the inverse image $f^{-1}(a)$?  
Since $\C$ is noncompact, there may be infinitely many points in the 
pre-image $f^{-1}(a)$.  Sometimes $f^{-1}(a)$ will be empty, meaning that 
$f$ misses the point $a$ in $\C P^1$.

The exponential map $\exp : \C \to \C P^1$ misses exactly two points, $0$ 
and $\infty$, in $\C P^1$.
According to a classical theorem of Picard, a nonconstant holomorphic map 
$f: \C \to \C P^1$ cannot miss more than two points.

Nevanlinna theory refines Picard's theorem in a beautiful way.  To each 
$a \in \C P^1$, it attaches a real number $\delta (a)$ between $0$ and $1$ 
inclusive, the \emph{deficiency index} of $a$. 
The deficiency index is a normalized way of 
counting the number of points in the inverse image.  If $f^{-1}(a)$ is 
empty, then the deficiency index is $1$.

In this context the first main theorem of Nevanlinna theory says that 
a nonconstant holomorphic map $f: \C \to \C P^1$ has 
deficiency index $0$ almost 
everywhere.  The second main theorem yields the stronger inequality:
\[
\sum_{a\in \C P^1} \delta (a) \le 2.
\]

Ahlfors generalized these two theorems to holomorphic maps with values in a 
complex projective space $\C P^n$.   

In [38] Bott and Chern souped up Nevanlinna's hard analysis to give a more 
conceptual proof of the first main theorem.

A by-product of Bott and Chern's excursion in Nevanlinna theory is the 
notion of a refined Chern class, now called the \emph{Bott-Chern class}, 
that has since been transformed into a 
powerful tool in Arakelov geometry and other aspects of modern number 
theory.

Briefly, the Bott-Chern classes arise as follows.  On a complex 
manifold $M$ 
the exterior derivative $d$ decomposes into a sum $d= \partial + 
\bar{\partial}$, and the smooth $k$-forms decompose into a direct 
sum of $(p,q)$-forms.  Let $A^{p,p}$ be the space of smooth $(p,p)$-forms 
on $M$.  Then the operator $\partial \bar{\partial}$ makes $\oplus A^{p,p}$ 
into a differential complex.  Thus, the cohomology $H^*\{ A^{p,p},
\partial \bar{\partial} \}$ is defined.

A Hermitian structure on a holomorphic rank $n$ vector 
bundle $E$ on $M$ determines a 
unique connection and hence a unique curvature tensor.  If $K$ and $K'$ are 
the curvature forms determined by two Hermitian structures on $E$ and 
$\phi$ is a $\GL (n,\C)$-invariant polynomial on $\gl (n,\C)$, 
then it is well known 
that $\phi (K)$ and $\phi(K')$ are global closed forms on $M$ 
whose difference is exact: 
\[
\phi (K) - \phi(K') = d \alpha
\]
for a differential form $\alpha$ on $M$.  This allows one to define the 
characteristic classes of $E$ as cohomology classes in $H^*(M)$.

In the holomorphic case, $\phi (K)$ and $\phi (K')$ are  
$(p,p)$-forms closed under $\partial \bar{\partial}$.  Bott and Chern 
found that in fact,
\[
\phi (K) - \phi(K') = \partial \bar{\partial} \beta
\]
for some $(p-1,p-1)$-form $\beta$.  
For a holomorphic vector bundle $E$, the Bott-Chern class of $E$ 
associated to an invariant polynomial $\phi$ is the cohomology class of 
$\phi (E)$, not in the usual 
cohomology, but in the cohomology of the complex $\{ A^{p,p},
\partial \bar{\partial}  \}$.

\section{Characteristic numbers and the Bott residue}

According to the celebrated Hopf index theorem, the Euler characteristic 
of a smooth manifold is equal to the number of zeros of a vector field on 
the manifold, each counted with its index.  In [41] and [43], Bott 
generalized the Hopf index theorem to other characteristic numbers such as 
the Pontryagin numbers of a real manifold and the Chern numbers of a 
complex manifold.

We will describe Bott's formula only for Chern numbers.  Let $M$ be a 
compact complex manifold of dimension $n$, and $c_1(M), \dots, c_n(M)$ the 
Chern classes of the tangent bundle of $M$.  The Chern numbers of $M$ are 
the integrals $\int_M \phi (c_1(M), \dots , c_n(M))$, as $\phi$ ranges 
over all weighted homogeneous polynomials of degree $n$.  Like the Hopf 
index theorem, Bott's formula computes a Chern number in terms of the zeros 
of a vector field $X$ on $M$, but the vector field must be holomorphic and the 
counting of the zeros is a little more subtle.

For any vector field $Y$ and any $C^{\infty}$ function $f$ on $M$, the Lie 
derivative $\liex$ satisfies:
\[
\liex (fY)= (Xf)Y+f\liex Y.
\]
It follows that at a zero $p$ of $X$,
\[
(\liex fY)_p = f(p) (\liex Y)p.
\]
Thus, at $p$, the Lie derivative $\liex$ induces an endomorphism
\[
L_p : T_pM \to T_p M
\]
of the tangent space of $M$ at $p$.  The zero $p$ is said to be 
\emph{nondegenerate} if $L_p$ is nonsingular.

For any endomorphism $A$ of a vector space $V$, we define the numbers 
$c_i (A)$ to be the coefficients of its characteristic polynomial:
\[
\det (I+tA) = \sum c_i (A) t^i.
\]

Bott's Chern number formula is as follows.  Let $M$ be a compact complex 
manifold of complex dimension $n$ and $X$ a holomorphic vector field having 
only isolated nondegenerate zeros on $M$.  For any weighted homogeneous  
polynomial $\phi (x_1, \dots, x_n)$, $\deg x_i = 2i$, 
\begin{equation} \label{characteristic}
\int_M \phi (c_1(M), \dots, c_n(M)) = \sum_p \dfrac{\phi (c_1 (L_p), \dots, 
c_n(L_p))}{c_n(L_p)},
\end{equation}
summed over all the zeros of the vector field. 
Note that by the definition of a nondegenerate zero, $c_n(L_p)$, which is 
$\det L_p$, is nonzero. 

In Bott's formula, if the polynomial $\phi$ does not have degree $2n$, then 
the left-hand side of \eqref{characteristic} is zero, and the formula 
gives an identity among the numbers $c_i(L_p)$.  For the polynomial $\phi 
(x_1, \dots, x_n) = x_n$, Bott's formula recovers the Hopf index theorem:
\[
\int_M c_n (M) = \sum_p \dfrac{c_n(L_p)}{c_n(L_p)} = \# \text{ zeros of } X.
\]

Bott's formula \eqref{characteristic} is reminiscent of Cauchy's residue 
formula and so the right-hand side of \eqref{characteristic} may be viewed 
as a residue of $\phi$ at $p$.

In [43] Bott generalized his Chern number formula \eqref{characteristic}, 
which assumes isolated zeros, to holomorphic vector fields with 
higher-dimensional zero sets and to bundles other than the tangent bundle 
(a vector field is a section of the tangent bundle). 

\section{The Atiyah-Bott fixed point theorem}

A continuous map of a finite polyhedron, $f:P \to P$, has a 
\emph{Lefschetz number}:
\[
L(f) = \sum (-1)^i \tr f^* | H^i (P),
\]
where $f^*$ is the induced homomorphism in cohomology and $\tr$ denotes the 
trace.
According to the Lefschetz fixed point theorem, if the Lefschetz number of 
$f$ is not zero, then $f$ has a fixed point.

In the smooth category the Lefschetz fixed point theorem has a quantitative
refinement.  A smooth map $f: M\to M$ from a compact manifold to itself is 
\emph{transversal} if its graph is transversal to the diagonal $\Delta$ in 
$M \times M$.  Analytically, $f$ is transversal if and only if at each 
fixed point $p$, 
\[
\det (1- f_{*,p}) \ne 0,
\]
where $f_{*,p}: T_pM \to T_pM$ is the differential of $f$ at $p$.

\begin{figure}[!h]
\begin{center}
\includegraphics[scale=0.3]{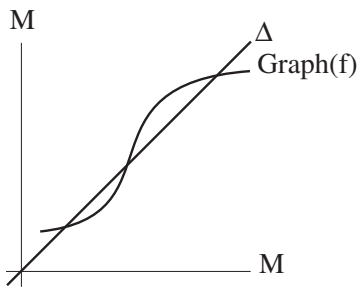}
 \caption{A transversal map $f$}
 \label{transversal}
\end{center}
\end{figure}

The $C^{\infty}$ Lefschetz fixed point theorem states that the Lefschetz 
number of a transversal map $f$ is the number of fixed points $f$ counted 
with multiplicity $\pm 1$ depending on the sign of the determinant
$\det (1- f_{*,p})$:
\[
L(f) = \sum_{f(p)=p} \pm 1.
\]

In the Sixties Atiyah and Bott proved a far-reaching generalization of the 
Lefschetz fixed point theorem ([42], [44]).  This type of result, relating a global 
invariant to a sum of local contributions, is a recurring theme in some of 
Bott's best work.

To explain it, recall that the real singular cohomology of $M$ is 
computable from the de Rham complex
\[
\Gamma (\Lambda^0) \overset{d}{\to} \Gamma (\Lambda^1) \overset{d}{\to}
\Gamma (\Lambda^2) \overset{d}{\to} \to \dots,
\]
where $\Lambda^q = \Lambda^q T^*M$ is the $q$th exterior power of the 
cotangent bundle.  The de Rham complex is an example of an
\emph{elliptic complex} on a manifold.

Let $E$ and $F$ be vector bundles of ranks $r_E$ and $r_F$ respectively 
over $M$.  An $\R$-linear map
\[
D: \Gamma (E) \to \Gamma (F)
\]
is a \emph{differential operator} if about every point in $M$
 there is a coordinate chart 
$(U, x_1, \dots, x_n)$ and trivializations for $E$ and $F$ relative to 
which $D$ can be written in the form
\[
D=\sum_{|\alpha| \le m} A^{\alpha}(x) \pd{x^{\ga}}, \quad x \in U, \quad
\quad \pd{x^{\alpha}} = \left(\pd{x_1}\right)^{\ga _1} \dots 
\left(\pd{x_n}\right)^{\ga _n} ,
\]
where $|\alpha|=\sum \alpha_i$ and
$A^{\ga}(x)$ is an $r_F \times r_E$ matrix that depends on $x$.  The 
order of $D$ is the highest $| \ga |$ that occurs.

Given a cotangent vector $\xi = \sum \xi_i dx_i \in T_x^* M$, we write
\[
\xi_{\ga} = \xi_1^{\ga _1} \dots \xi_n^{\ga _n}
\]
and define the \emph{symbol} of a differential operator $D$ of order $m$ 
to be
\[
\sigma (D, \xi)_x = \sum_{|\alpha| = m} A^{\ga} (x) \xi _{\ga} \in \Hom 
(E_x, F_x).
\]
In other words, the symbol of $D$ is obtained by first discarding all but 
the highest-order terms of $D$ and then replacing $\partial / \partial 
x^{\ga}$ by $\xi_{\ga}$.  Because $\xi_i$ transforms like $\partial / \partial 
x^i$ under a change of coordinates, it is not difficult to show that 
the symbol is well-defined, independent of the coordinate system.

Let $E_i$ be vector bundles over a manifold $M$.  A differential complex
\begin{equation} \label{complex}
\mathcal{E}: 0 \to \Gamma(E_0) \overset{D}{\to} \Gamma(E_1) \overset{D}{\to}
\Gamma(E_0) \overset{D}{\to}  \dots, \quad D^2=0,
\end{equation}
is \emph{elliptic} if for each nonzero cotangent vector $\xi \in T_x^*M$, 
the associated symbol sequence
\[
0 \to E_{0,x} \xrightarrow{\sigma(D,\xi)} E_{1,x} \xrightarrow{\sigma(D,\xi)}
E_{2,x} \xrightarrow{\sigma(D,\xi)} \dots
\]
is an exact sequence of vector spaces.

A fundamental consequence of ellipticity is that all the cohomology spaces 
$H^i = H^i(\Gamma (E_{*}))$ are finite-dimensional.

An \emph{endomorphism} of the complex \eqref{complex} is a collection of linear 
maps $T_i: \Gamma (E_i) \to \Gamma (E_i)$ such that
\[
T_{i+1} \comp D = D \comp T_i
\]
for all $i$.  Such a collection $T=\{T_i\}$ induces maps in cohomology
$T_i^*:H^i \to H^i$.  The \emph{Lefschetz number} of $T$ is then defined to be
\[
L(T) = \sum (-1)^i \tr T_i^*.
\]

A map $f: M \to M$ induces a natural map
\[
\Gamma_f : \Gamma (E) \to \Gamma(f\inv E)
\]
by composition: $\Gamma_f (s) = s \comp f$.  There is no natural way to 
induce a map of sections: $\Gamma (E) \to \Gamma (E)$.
However, if there is a bundle map $\phi: f\inv E \to E$, then the composite
\[
\Gamma (E) \overset{\Gamma_f}{\to} \Gamma(f\inv E) 
\overset{\tilde{\phi}}{\to} \Gamma (E)
\]
is an endomorphism of $\Gamma (E)$.  Any bundle map $\phi : f\inv E \to E$ 
is called a \emph{lifting} of $f$ to $E$.  At each point $x\in M$, a 
lifting $\phi$ is nothing other than a linear map $\phi_x: E_{f(x)} \to 
E_x$.

In the case of the de Rham complex, a map $f:M\to M$ induces a linear map
$f_x^*: T_{f(x)}^*M \to T_x^*M$ and hence a linear map
\[
\Lambda^q f_x^* : \Lambda^q T_{f(x)}^*M \to \Lambda T_x^* M,
\]
which is the lifting that finally defines the pullback of differential 
forms $f^* : \Gamma(\Lambda^q T^* M) \to \Gamma (\Lambda^q T^* M)$.

\begin{theorem}[Atiyah-Bott fixed point theorem]
Given an elliptic complex \eqref{complex} on a compact manifold $M$, 
suppose $f: M \to M$ has a lifting $\phi_i : f\inv E_i \to E_i$ for each 
$i$ such that the induced maps $T_i: \Gamma(E_i) \to \Gamma (E_i)$ give an 
endomorphism of the elliptic complex.  Then the Lefschetz number of $T$ is
given by
\[
L(T)= \sum_{f(x)=x} \dfrac{\sum (-1)^i \tr \phi_{i,x}}
{|\det(1- f_{*,x})|}.
\]
\end{theorem}

As evidence of its centrality, the Atiyah-Bott fixed point theorem has an 
astonishing range of applicability.

Here is an easily stated corollary in algebraic geometry:  any holomorphic 
map of a rational algebraic manifold to itself has a fixed point.

Specializing the Atiyah-Bott fixed point theorem to the de Rham complex, 
one recovers the classical Lefschetz fixed point theorem.  When applied 
to other geometrically interesting elliptic complexes, Atiyah and Bott 
obtained new fixed point theorems, such as a holomorphic Lefschetz fixed 
point theorem in the complex analytic case and a signature formula in the 
Riemannian case.  In the homogeneous case, the fixed point theorem implies 
the Weyl character formula.

\section{Obstruction to integrability}

A subbundle $E$ of the tangent bundle $TM$ of a manifold $M$ assigns to 
each point $x$ of the manifold a subspace $E_x$ of the tangent space 
$T_xM$.  An \emph{integrable manifold} of the subbundle $E$ is a 
submanifold $N$ of $M$ whose tangent space $T_xN$ at each point $x$ 
in $N$ is $E_x$.  The subbundle $E$ is said to be \emph{integrable} if 
for each point $x$ in $M$, there is an integrable manifold of $E$ passing 
through $x$.

By the Frobenius theorem, often proven in a first-year graduate course, a 
subbundle $E$ of the tangent bundle $TM$ is integrable if and only if its 
space of sections $\Gamma (E)$ is closed under the Lie bracket.

The Pontryagin ring $\Pont(V)$ 
of a vector bundle $V$ over $M$ is defined to be the 
subring of the cohomology ring $H^*(M)$ generated by the Pontryagin 
classes of the bundle $V$.  In [51] Bott found an obstruction to the 
integrability of $E$ in terms of the Pontryagin ring of the quotient bundle
$Q:=TM/E$.  More precisely, if a subbundle $E$ of the tangent 
bundle $TM$ is integrable, then the Pontryagin ring $\Pont (Q)$ vanishes in 
dimensions greater than twice the rank of $Q$.

What is so striking about this theorem is not only the simplicity of the 
statement, but also the simplicity of its proof.  It spawned tremendous 
developments in foliation theory in the Seventies, as recounted in 
\cite{conlon} and \cite{haefliger95}.

\section{The cohomology of the vector fields on a manifold}

For a finite-dimensional Lie algebra $L$, let $A^q(L)$ be the space of 
alternating $q$-forms on $L$.  Taking cues from the Lie algebra of 
left-invariant vector fields on a Lie group, one defines the 
differential
\[
d: A^q(L) \to A^{q+1}(L)
\]
by
\begin{equation} \label{differential}
(d\omega)(X_0, \dots, X_q) = \sum_{i < j} (-1)^{i+j} \omega([X_i, X_j], 
X_0, \dots, \hat{X_i}, \dots, \hat{X_j}, \dots, X_q).
\end{equation}
As usual, the hat $\hat{}$ over $X_i$ means that $X_i$ is to be omitted.
This makes $A^*(L)$ into a differential complex, whose cohomology is by 
definition the cohomology of the Lie algebra $L$. 

When $L$ is the infinite-dimensional Lie algebra $L(M)$ of vector fields 
on a manifold $M$, the formula \eqref{differential} still makes sense, 
but the space of all alternating forms $A^*(L(M))$ is too large for its 
cohomology to be computable.  Gelfand and Fuks proposed 
putting a topology, the $C^{\infty}$ 
topology, on $L(M)$, and computing instead the cohomology of 
the \emph{continuous} alternating 
forms on $L(M)$.  The \emph{Gelfand-Fuks cohomology} of $M$ is the 
cohomology of the complex $\{ A_c^*(L(M)), d\}$ of continuous forms.  They 
hoped to find in this way new invariants of a manifold.  As an example, they 
computed the Gelfand-Fuks cohomology of a circle.

It is not clear from the definition that the Gelfand-Fuks cohomology is a 
homotopy invariant.  In [71] Bott and Segal proved that the Gelfand-Fuks 
cohomology of a manifold $M$ is the singular cohomology of a space 
functorially constructed from $M$.  Haefliger \cite{haefliger76} and Trauber gave a very 
different proof of this same result.  The homotopy invariance of the 
Gelfand-Fuks cohomology follows.  At the same time it also showed that the 
Gelfand-Fuks cohomology produces no new invariants.

\section{Localization in equivariant cohomology}

Just as singular cohomology is a functor from the category of topological 
spaces to the category of rings, so when a group $G$ acts on a space $M$, 
one seeks a functor that would incorporate both the topology of the space 
and the action of the group.

The naive construction of taking the cohomology of the quotient space 
$M/G$ is unsatisfactory because for a nonfree action, the topology of the 
quotient can be quite bad.  A solution is to find a contractible space 
$EG$ on which $G$ acts freely, for then $EG \times M$ will have the same 
homotopy type as $M$ and the group $G$ will act freely on $EG \times M$ via 
the diagonal action.  It is well known that such a space is the total 
space of the universal $G$-bundle $EG\to BG$, whose base space is the 
classifying space of $G$.  
The homotopy theorists have defined the \emph{homotopy quotient} 
$M_G$ of $M$ by $G$ to be the quotient 
space $(EG \times M)/G$, and the \emph{equivariant cohomology} $H_G^*(M)$ 
to be the ordinary cohomology of its homotopy quotient $M_G$.

The equivariant cohomology of the simplest $G$-space, a point, is already 
quite interesting, for it is the ordinary cohomology of the classifying 
space of $G$:
\[
H_G^*(\pt)= H^*((EG\times \pt)/G )= H^*(EG/G)= H^*(BG).
\]

Since equivariant cohomology is a functor of $G$-spaces, the constant map 
$M \to \pt$ induces a homomorphism $H_G^*(\pt) \to H_G^*(M)$.  Thus, the 
equivariant cohomology $H_G^*(M)$ has the structure of a module over 
$H^*(BG)$. 

Characteristic classes of vector bundles over $M$ extend to equivariant 
characteristic classes of equivariant vector bundles.

When $M$ is a manifold, there is a push-forward map $\pi_*^M: 
H_G^*(M) \to H_G^*(\pt)$, akin to integration along the fiber.

Suppose a torus $T$ acts 
on a compact manifold $M$ with fixed point set $F$, and $\phi \in H_T^*(M)$ 
is an equivariantly closed class.  Let $P$ be the connected components of 
$F$ and let $\iota_P:P\to M$ be the inclusion map, $\nu_P$ the normal 
bundle of $P$ in $M$, and $e(\nu_P)$ the equivariant Euler class of 
$\nu_P$.  In [82] Atiyah and Bott proved a localization theorem for the 
equivariant 
cohomology $H_T^*(M)$ with real coefficients: 
\[
\pi_*^M \phi = \sum_P \pi_*^P \left( \dfrac{\iota_P^* \phi}{e(\nu_P)} 
\right).
\]
It should be noted that Berline and Vergne \cite{berline-vergne} 
independently proved the same theorem at about the same time. 

This localization theorem has as consequences the following results of 
Duistermaat and Heckman on a symplectic manifold $(M,\omega)$ of dimension 
of $2n$:
\begin{enumerate}
\item[1)] If a torus action on $M$ preserves the symplectic form and has a 
moment map $f$, then the push-forward $f_{*} (\omega^n)$ 
of the symplectic volume under the 
moment map is piecewise polynomial.
\item[2)] Under the same hypotheses, the stationary phase approximation for
the integral
\[
\int_M e^{-itf} \dfrac{\omega^n}{n!}
\]
is exact.
\end{enumerate}

In case the vector field on the manifold is generated 
by a circle action, the localization theorem
 specializes to Bott's Chern number formulas 
[41] of the Sixties,  thus providing an alternative explanation for the Chern 
number formulas.

\section{Yang-Mills equations over Riemann surfaces}

In algebraic geometry it is well known that for any degree $d$ the set of 
isomorphism classes of holomorphic line bundles of degree $d$ over a 
Riemann surface $M$ of genus $g$ forms a smooth projective variety which is 
topologically a torus of dimension $g$.  This space is called the 
\emph{moduli space} of holomorphic line bundles of degree $d$ over $M$.

For holomorphic vector bundles of rank $k \ge 2$, the situation is far 
more complicated.  First, in order to have an algebraic structure on the 
moduli space, it is necessary to discard the so-called ``unstable'' bundles 
in the sense of Mumford.  It is then known that for $k$ and $d$ relatively 
prime, the isomorphism classes of the remaining bundles, called 
``semistable bundles,'' form a smooth projective variety $N(k,d)$.

In \cite{newstead} Newstead computed the Poincar\'e polynomial of 
$N(2,1)$.  Apart from this, the topology of $N(k,d)$ 
remained mysterious. 

In [81] Atiyah and Bott introduced the new and powerful method of 
equivariant Morse theory to study the topology of these moduli spaces.

Let $P= M \times U(n)$ be the trivial principal $U(n)$-bundle over the 
Riemann surface $M$, $\cala=\cala(P)$ the affine space of connections on 
$P$, and $\calg=\calg(P)$ the gauge group, i.e.,
 the group of automorphisms of $P$ that cover the 
identity.  Then the gauge group $\calg (P)$ acts on the space 
$\cala (P)$ of connections and there is a Yang-Mills 
functional $L$ on $\cala(P)$ invariant under the action of the gauge group. 

Equivariant Morse theory harks back to Bott's extension of classical Morse 
theory to nondegenerate critical manifolds three decades earlier.
The key result of Atiyah and Bott is that the Yang-Mills functional $L$ is a 
perfect equivariant Morse function on $\cala (P)$.  This means the 
equivariant Poincar\'e series of $\cala (P)$ is equal to the 
equivariant Morse series of $L$:
\begin{equation} \label{eqmorse}
P_t^{\calg} (\cala(P)) = \calm_t^{\calg}(L).
\end{equation}

Once one unravels the definition, the left-hand side of \eqref{eqmorse} is 
simply the Poincar\'e series of the classifying space of $\calg (P)$, which 
is computable from homotopy considerations.  The right-hand side of 
\eqref{eqmorse} is the sum of contributions from all the critical sets of 
$L$.  By the work of Narasimhan and Seshadri, the minimum of $L$ is 
precisely the moduli space $N(k,d)$.  It contributes its Poincar\'e 
polynomial to the equivariant Morse series of $L$.  By an inductive 
procedure, Atiyah and Bott were able to compute the contributions of all 
the higher critical sets.  They then solved \eqref{eqmorse} for the 
Poincar\'e polynomial of $N(k,d)$.

\section{Witten's rigidity theorem}

Let $E$ and $F$ be vector bundles over a compact manifold $M$.  If a 
differential operator $D: \Gamma (E) \to \Gamma (F)$ is elliptic, then 
$\ker D$ and $\coker D$ are finite-dimensional vector spaces and we can 
define the \emph{index} of $D$ to be the virtual vector space
\[
\index D = \ker D - \coker D.
\]

Now suppose a Lie group $G$ acts on $M$, and $E$ and $F$ are 
$G$-equivariant vector bundles over $M$.  Then $G$ acts on $\Gamma (E)$ by
\[
(g.s)(x)= g.(s(g\inv .x)),
\]
for $g \in G$, $s\in \Gamma(E), x\in M$.  The $G$-action is said to 
\emph{preserve} the differential operator $D$ if the actions of $G$ on 
$\Gamma (E)$ and $\Gamma (F)$ commute with $D$.
In this case $\ker D$ and $\coker D$ are representations of $G$, and so 
$\index D$ is a virtual representation of $G$.  We say that the operator 
$D$ is \emph{rigid} if its index is a multiple of the trivial 
representation of dimension $1$.  The rigidity of $D$ means that any 
nontrivial irreducible representation of $G$ in $\ker D$ occurs in $\coker 
D$ with the same multiplicity and vice versa.

If the multiple $m$ is positive, then $m.1= 1\oplus  \dots \oplus 1$ is 
the trivial representation of dimension of $m$.  If $m$ is negative, the 
$m.1$ is a virtual representation and the rigidity of $D$ implies that the 
trivial representation $1$ occurs more often in $\coker D$ than in $\ker D$.

For a circle action on a compact oriented Riemannian manifold, it is well 
known that the Hodge operator $d+d*: \Omega^{\text{even}} \to 
\Omega^{\text{odd}}$ and the signature operator $d_s=d+d^*:\Omega^+ \to 
\Omega^-$ are both rigid.

An oriented Riemannian manifold of dimension $n$ has an atlas whose transition 
functions take values in $\SO (n)$.  The manifold is called a \emph{spin 
manifold} if it is possible to lift the transition functions to the double cover 
$\Spin (n)$ of $\SO (n)$.

Inspired by physics, Witten discovered infinitely many rigid elliptic 
operators on a compact spin manifold with a circle action.  
They are typically of the form 
$d_s \otimes R$, where $d_s$ is the signature operator
and $R$ is some combination of the exterior and the symmetric powers of 
the tangent bundle.  In [91] Bott and Taubes  
found a proof, more accessible to mathematicians, of Witten's results, by 
recasting the rigidity theorem as a consequence of 
the Atiyah-Bott fixed point theorem.

The idea of [91] is as follows.  To decompose a representation, one needs 
to know only its trace, since the trace determines the representation.
By assumption, the action of $G$ on the elliptic complex $D: \Gamma (E) \to 
\Gamma (F)$ commutes with $D$.  This means each element $g$ in $G$ is an 
endomorphism of the elliptic complex.  It therefore induces an 
endomorphism $g^*$ in the cohomology of the complex.  But $H^0 = \ker D$ and 
$H^1=\coker D$.  The alternating sum of the trace of $g^*$ in cohomology 
is precisely the left-hand side of the Atiyah-Bott fixed point theorem.  
It then stands to reason that the fixed point theorem could be used to 
decompose the index of $D$ into irreducible representations.

\bigskip
\noindent
{\bf Papers of Raoul Bott discussed in this article}

\medskip
\begin{enumerate}
\item[{[1]}]  (with R. J. Duffin) Impedance synthesis without use of 
transformers, J.\ Appl.\ Phys., {\bf 20} (1949), 816.
\item[{[8]}] On torsion in Lie groups, Proc.\ NAS, {\bf 40} (1954), 586--588.
\item[{[9]}] Nondegenerate critical manifolds, Ann.\ of Math.\ {\bf 60} 
(1954), 248--261.
\item[{[12]}] (with H. Samelson) The cohomology ring of $G/T$, 
Proc.\ NAS, {\bf 41} (1955), 490--493.
\item[{[14]}] On the iteration of closed geodesics and the Sturm 
Intersection theory, Comm.\ Pure Appl.\ Math.\ {\bf IX} (1956), 171--206.
\item[{[15]}] Homogeneous vector bundles, Ann.\ of Math.\ {\bf 66} (1957), 
933--935.
\item[{[24]}] The stable homotopy of the classical groups, Ann.\ of Math.\ 
{\bf 70} (1959), 313--337.
\item[{[33]}] (with M.~F.~Atiyah and A.~Shapiro) Clifford modules, 
Topology {\bf 3} (1965), 3--38.
\item[{[37]}] The index theorem for homogeneous differential operators, in: 
{\it Differential and Combinatorial Topology: A Symposium in Honor of 
Marston Morse}, Princeton, (1964), 167--186.
\item[{[38]}] (with S.~Chern) Hermitian vector bundles and the 
equidistribution of the zeroes of their holomorphic sections, Acta 
Mathematics {\bf 114} (1964), 71--112.
\item[{[41]}] Vector fields and characteristic numbers, Mich.\ Math.\ 
Jour. {\bf 14} (1967), 231--244.
\item[{[42]}] (with M.~F.~Atiyah) A Lefschetz fixed point formula for 
elliptic complexes: I, Ann.\ of Math.\ {\bf 86} (1967), 374--407.
\item[{[43]}] A residue formula for holomorphic vector fields, J.\ 
Differential Geom.\ {\bf 1} (1967), 311--330.
\item[{[44]}] (with M.~F.~Atiyah) A Lefschetz fixed point formula for 
elliptic complexes: II, Ann.\ of Math.\ {\bf 88} (1968), 451--491.
\item[{[51]}] On a topological obstruction to integrability, in: {\it 
Global Analysis}, Proceedings of Symposia in Pure Math.\ {\bf XVI} (1970), 
127--131.
\item[{[71]}] (with G. Segal) The cohomology of the vector fields on a 
manifold, Topology {\bf 16} (1977), 285--298.
\item[{[81]}] (with M.~F.~Atiyah) The Yang-Mills equations over Riemann 
surfaces, Phil.\ Trans.\ R. Soc.\ Lond.\ {\bf A 308} (1982), 524--615.
\item[{[82]}] (with M.~F.~Atiyah) The moment map and equivariant 
cohomology, Topology {\bf 23} (1984), 1--28.
\item[{[91]}] (with C. Taubes) On the rigidity theorem of Witten, 
Jour.\ of the Amer.\ Math.\ 
Soc.\ {\bf 2} (1989), 137--186. 
\end{enumerate} 

\bigskip
\noindent
{\bf Update to the Bibliography of Raoul Bott}

\medskip
\noindent
Raoul Bott's Bibliography in his Collected Papers \cite{collected} 
is not complete.
This update completes the Bibliography as of June 2000.

\medskip
\begin{enumerate}
\item[{[20]}] In memoriam Sumner B.~Myers, Mich.\ Math.\ Jour., {\bf 5} (1958), 
1--4.
\item[{[79]}] (with L. Tu) {\it Differential Forms in Algebraic Topology}, 
Springer-Verlag, (1982), 1--331.
\item[{[94]}] Georges de Rham: 1901--1990, {\it Notices Amer.\ Math.\ Soc.} 
{\bf 38} (1991), 114--115.
\item[{[95]}] Stable bundles revisited, Surveys in Differential Geometry, 
{\bf 1} (1991), 1--18.
\item[{[96]}] On E.~Verlinde's formula in the context of stable bundles, 
Internat. J. Modern Phys. A {\bf 6} (1991), 2847--2858.
\item[{[97]}] Nomination for Stephen Smale, {\it Notices Amer.\ Math.\ Soc.},
{\bf 38} (1991), 758--760..
\item[{[98]}] On knot and manifold invariants, in:  {\it New Symmetry 
Principles in Quantum Field Theory (Carg\`ese, 1991)}, NATO Adv.\ Sci.\ 
Inst.\ Ser.\ B Phys.\ 295, Plenum, (1992), 
37--52.
\item[{[99]}] Topological aspects of loop groups, in: {\it Topological 
Quantum Field Theories and Geometry of Loop Spaces (Budapest, 1989)}, 
World Sci.\ Publishing, (1992), 65--80.
\item[{[100]}] For the Chern volume, in: {\it Chern---a Great Geometer of the 
Twentieth Century}, Internat.\ Press, (1992), 106--108.
\item[{[101]}] Reflections on the theme of the poster, in: {\it Topological 
Methods in Modern Mathematics (Stony Brook, NY, 1991)}, Publish or Perish 
(1993), 125--135.
\item[{[102]}] Luncheon talk and nomination for Stephen Smale, in: {\it 
From Topology to Computation: Proceedings of the Smalefest (Berkeley, CA, 
1990)}, Springer, (1993), 67--72.
\item[{[103]}] (with C. Taubes) On the self-linking of knots, J. Math.\ 
Phys.\ {\bf 35} (1994), 5247--5287.
\item[{[104]}] On invariants of manifolds, in: {\it Modern Methods in 
Complex Analysis (Princeton, NJ, 1992)}, Ann.\ of Math.\ Stud., 137, 
Princeton Univ. Press, (1995), 29--39.
\item[{[105]}] Configuration spaces and imbedding invariants, Turkish J.\ 
Math.\ {\bf 20} (1996), 1--17.
\item[{[106]}] Configuration spaces and imbedding problems, in: {\it 
Geometry and Physics (Aarhus, 1995)}, Lecture Notes in Pure and Appl. 
Math., {\bf 184}, Dekker (1997), 135--140.
\item[{[107]}] Critical point theory in mathematics and in mathematical 
physics, Turkish J.\ Math.\ {\bf 21} (1997), 9--40.
\item[{[108]}] (with A. Cattaneo) Integral invariants of $3$-manifolds, J.\ 
Differential Geom.\ {\bf 48} (1998), 91--133.
\item[{[109]}] An introduction to equivariant cohomology, in: {\it Quantum 
Field Theory: Perspective and Prospective}, Kluwer Academic Publishers, 
(1999), 35--57.
\item[{[110]}] (with A. Cattaneo) Integral invariants of $3$-manifolds, 
II, to appear in  J.\ Differential Geom.
\end{enumerate}


\begin{thebibliography}{B9}

\bibitem[BV]{berline-vergne} N. Berline and M. Vergne, 
Classes caract\'eristiques \'equivariantes, 
Formule de localisation en cohomologie \'equivariante, 
C. R. Acad.\ Sci.\ Paris S\'er.\ I Math.\ {\bf 295} (1982), No.\ 9, 539--541.

\bibitem[B1]{recoll} R. Bott, Some recollections from 30 years ago, in: {\it 
Constructive Approaches to Mathematical Models}, Academic Press, 1979, 
33--39.

\bibitem[B2]{specker} R. Bott, An equivariant setting of the Morse theory, 
L'Enseignement math\'ematique, {\bf XXVI} (1980), 68--75.

\bibitem[B3]{dioszeger} R. Bott, The Dioszeger Years (1923--1929), 
in: {\it Raoul Bott: 
Collected Papers}, Vol.\ 1, Birkh\"auser, Boston, 1994, 11--26.

\bibitem[B4]{auto} R. Bott, Autobiographical Sketch, in {\it Raoul Bott: 
Collected Papers}, Vol.\ 1, Birkh\"auser, Boston, 1994, 3--9.

\bibitem[B5]{collected} R. Bott, {\it Collected Papers}, Vol.\ 1--4, 
Birkh\"auser, Boston, 
1994, 1995.

\bibitem[C]{conlon} L. Conlon, Raoul Bott, foliations, and characteristic 
classes: an appreciation, in {\it Raoul Bott: Collected Papers}, Vol.\ 3, 
Birkh\"auser, Boston, 1995, xxiv--xxvi.

\bibitem[H]{haefliger76} A. Haefliger, Sur la cohomologie de Gel'fand-Fuks, 
Ann.\ scient.\ Ec.\ Norm.\ Sup.\ {\bf 9} (1976), 503--532.

\bibitem[H1]{haefliger95} A. Haefliger, Raoul Bott and foliation theory in 
the 1970s, in {\it Raoul Bott: Collected Papers}, Vol.\ 3, 
Birkh\"auser, Boston, 1995, xxvii--xxxi. 

\bibitem[N]{newstead} P. E. Newstead, Characteristic classes of stable 
bundles over an algebraic curve, Trans.\ Am.\ Math.\ Soc.\ {\bf 169}, 337--345.  

\end{thebibliography}
\end{document}